\documentclass[10pt]{amsart}

\usepackage{amsfonts,amssymb}

\newtheorem{proposition}{Proposition}[section]
\newtheorem{lemma}[proposition]{Lemma}
\newtheorem{corollary}[proposition]{Corollary}
\newtheorem{theorem}[proposition]{Theorem}
\theoremstyle{definition}
\newtheorem{definition}[proposition]{Definition}
\newtheorem{example}[proposition]{Example}

\theoremstyle{remark}
\newtheorem{remark}[proposition]{Remark}

\newcommand{\thlabel}[1]{\label{th:#1}}

\newcommand{\thref}[1]{Theorem~\ref{th:#1}}

\newcommand{\lelabel}[1]{\label{le:#1}}
\newcommand{\leref}[1]{Lemma~\ref{le:#1}}
\newcommand{\prlabel}[1]{\label{pr:#1}}
\newcommand{\prref}[1]{Proposition~\ref{pr:#1}}
\newcommand{\colabel}[1]{\label{co:#1}}
\newcommand{\coref}[1]{Corollary~\ref{co:#1}}

\newcommand{\relabel}[1]{\label{re:#1}}
\newcommand{\reref}[1]{Remark~\ref{re:#1}}

\newcommand{\delabel}[1]{\label{de:#1}}

\def\ra{\rightarrow}

\newcommand\smi{\mbox{$S^{-1}$}}
\def\ot{\otimes}
\def\va{\varepsilon}

\def\le{\langle}
\def\ri{\rangle}
\def\l{\lambda}

\def\va{\varepsilon}

\def\rh{\rightharpoonup}
\def\lh{\leftharpoonup}

\def\ra{\rightarrow}
\def\a{\alpha}

\newcommand{\sinv}{S^{-1}}

\begin{document}
\title[Radford's $S^4$ formula]
{Radford's $S^4$ formula for co-Frobenius Hopf algebras}
\author{Margaret Beattie}
\address{Department of Mathematics and Computer Science, Mount Allison University,
Sack-ville, NB E4L 1E6, Canada} \email{mbeattie@mta.ca}\thanks{The
first author's research is supported by NSERC}
\author{Daniel Bulacu}
\address{Faculty of Mathematics and Informatics, University
of Bucharest, Str. Academiei 14, RO-010014 Bucharest 1, Romania}
\email{dbulacu@al.math.unibuc.ro}
\thanks{
The second author had support from LIEGRITS, a Marie Curie
Research Training Network funded by the European community as
project MRTN-CT 2003-505078 and a PDF at Mount Allison University.
He thanks the University of Almeria (Spain), and Mount Allison
University (Canada)  for their warm hospitality.}
\author{Blas Torrecillas}
\address{Department of Algebra and Analysis,
University of Almeria, 04071 Almeria, Spain}
\email{btorreci@ual.es}
\thanks{The work of the third author was supported by grant MTM2005-3227 from
MEC} 
\begin{abstract}
This note extends Radford's formula for the fourth power of the
antipode of a finite dimensional Hopf algebra to co-Frobenius Hopf
algebras and studies equivalent conditions to a Hopf algebra being
involutory for   finite dimensional and co-Frobenius Hopf algebras.
\end{abstract}
\maketitle
%%%%%%%%%%%%%%%%%%%%%%%%%%%%%%%%%%%%%%%%%%%%%%
\section{Introduction and Preliminaries}
%%%%%%%%%%%%%%%%%%%%%%%%%%%%%%%%%%%%%%%%%%%%%%
${\;\;\;}$ A key result in the theory of finite dimensional Hopf
algebras is Radford's $S^4$ formula. It asserts that for $S$
  the antipode   of a finite dimensional Hopf
algebra $H$ over a field $k$, then  for all $h\in H$,
\begin{equation}\label{rsf}
S^4(h)=g(\a \rh h\lh \a ^{-1})g^{-1},
\end{equation}
where $g$ and $\a$ are the distinguished grouplike elements of $H$
and $H^*$, respectively. This formula was initially proved by
Larson in \cite[Theorem 5.5]{larson} for finite dimensional
unimodular Hopf algebras and   extended by Radford in
\cite[Proposition 6]{rad1} to any   finite dimensional Hopf
algebra $H$. The techniques used by Larson and Radford are
similar, having their roots in a paper of Sweedler \cite{sw2}.
Their idea was to compute the square, and then the fourth power,
of the antipode   in terms of a nonsingular bilinear form
associated  to $H$. More recently, Schneider in \cite[Theorem
3.8]{sch} used the fact that a finite dimensional Hopf algebra is
a Frobenius algebra to provide a more transparent proof. The dual
bases of $H$ in terms of integrals for $H$ and $H^*$ were
explicitly computed in two different ways. Then Radford's $S^4$
formula follows from the two expressions for the Nakayama
automorphism of $H$. Kadison and Stolin in \cite{ks} use this same
point of view  to obtain an analogous result for
  Hopf algebras over commutative rings that are
Frobenius algebras as do Doi and Takeuchi \cite{doitakeuchi} in
their proof of the $S^4$ formula for biFrobenius algebras.
Similarly, Montgomery \cite{moconstanta} uses this approach in her
discussion of finite-dimensional Hopf algebras and new proof that in
characteristic 0, a finite-dimensional Hopf algebra $H$ is
involutory if and only if $H$ or $H^*$ is semisimple.

Moreover,  the $S^4$ formula  was recently generalized for
quasi-Hopf algebras in \cite[Corollary 6.3]{hn3}, for weak Hopf
algebras in \cite[Theorem 5.13]{n} and for $QFH$-algebras
(quasi-Hopf algebras over commutative rings that are Frobenius
algebras) in \cite[Theorem 3.3]{ka}. A braided version of the $S^4$
formula is proved \cite{bklt} and in \cite{doitakeuchi}. Perhaps the
most striking generalization is the categorical $S^4$ formula
associated to any finite tensor category, see \cite[Theorem 3.3]{eno}.\\
${\;\;\;}$  We begin with a review of the spaces of integrals in
Hopf algebras and give a short proof of the $S^4$ formula for
finite dimensional Hopf algebras using integrals. We then prove an
analogue of the $S^4$ formula for co-Frobenius Hopf algebras,
providing another new proof for the finite dimensional case.  In
the last section, we study
conditions on left or right integrals which are equivalent to the Hopf
algebra being involutory.    \\
${\;\;\;}$ Key to our arguments   is the use of the one
dimensional spaces of  integrals in $H$ and   $H^*$ and versions
of the isomorphism induced by the structure theorem for Hopf
modules \cite[Theorem 4.1.1]{sw}. We   recall some basics about
spaces of integrals.

\subsection{Preliminaries: Spaces of Integrals} Throughout,
 $H$ will be either a finite dimensional or a co-Frobenius
Hopf algebra with antipode $S$. In either case,  $S$ is bijective
with composition inverse  denoted by $\sinv$. We will use the
Heyneman-Sweedler notation \cite{sw}, $\Delta (h)=h_1\ot h_2$, $h\in
H$, for comultiplication (summation understood). We assume that the
reader is familiar with the basic theory of Hopf algebras; see
\cite{dnr, mo, sw} for example. We work over a commutative field $k$
throughout.

Any Hopf algebra $H$ is an $H^*$-bimodule where for   $h^*\in H^*$
and $h\in H$,
\begin{equation}\label{H*action}
h^*\rh h=h^*(h_2)h_1,\mbox{  and  } h\lh h^*=h^*(h_1)h_2.
\end{equation}
 Also the algebra $H^*$ is an $H$-bimodule, where for $h^*\in H^*$ and $h, h'\in
 H$,
\begin{equation}\label{Haction}
\le h\rh h^*, h'\ri=\le h^*, h'h\ri ,\mbox{  and  } \le h^*\lh h,
h'\ri=\le h^*, hh'\ri.
\end{equation}
The fact that $S$ is an anti-algebra morphism yields that for $h\in
H$, $h^* \in H^*$,
\begin{equation}\label{Saction}
(h\rh h^*)\circ S=(h^* \circ S)\lh \smi (h),~~{\rm and}~~ (h^*\lh
h)\circ S=\smi (h)\rh (h^*\circ S).\end{equation}

Now suppose that $H$ is finite dimensional and recall the
definition of an integral in $H$. An element $t\in H$ is called a
left integral in $H$ if $ht=\va(h)t$ for all $h\in H$. The space
of left integrals is denoted $\int_l^H$ and is a one dimensional
ideal of $H$. Similarly if $T\in H$ has the property that $Th =
\epsilon(h)T$ for all $h \in H$, then $T$ is called a right
integral in $H$, and we write $T \in \int_r^H$.

Since $\int_l^H$ is a one dimensional ideal in $H$, there exists a
unique $\a \in H^*$ such that for any $t \in \int_l^H$ and for all
$h\in H$
\begin{equation}\label{alpha}
th=\a (h)t.
\end{equation}
 The map $\a $ is   grouplike in   $H^*$, the dual Hopf
 algebra of $H$, i.e.,
  $\alpha$ is an invertible algebra map from $H$ to $k$.   Following Radford
\cite{rad3}, we call  $\a $ the distinguished grouplike element of
$H^*$. Clearly, $\int_r^H = \int_l^H$ if and only if $\a =\va$. An
integral is called cocommutative if it is a cocommutative element of
$H$.

 If $0 \neq l$ is  a left or a right integral in $H$, then the
 maps
  from $H^*$ to $H$ given by
\begin{equation}\label{bij}
  h^* \mapsto (h^*\rh l=h^*(l_2)l_1) \mbox{   and   }
h^*\mapsto (l\lh h^* =h^*(l_1)l_2 )
\end{equation}
are bijections (see   \cite{dnr, sw}).

Now   let $H$ be a Hopf algebra which is not necessarily finite
dimensional. An element $\Lambda \in H^*$ is called a right
integral for $H$ in $H^*$ if $\Lambda(h_1)h_2 = \Lambda(h)1$ for
all $h \in H$.  The space of right integrals for $H$ in $H^*$ is
denoted $\int_r^{H^*}$. Similarly, $\lambda \in H^*$ is called a
left integral for $H$ in $H^*$ if $h_1\lambda(h_2) = \lambda(h)1$
for all $h \in H$; the space of left integrals is denoted
$\int_l^{H^*}$. The dimensions of the spaces $\int_l^{H^*}$ and
$\int_r^{H^*}$ are equal and at most 1.  For $H$   finite
dimensional, these are just the one dimensional spaces of
integrals for   $H^*$. An
  integral $\Gamma$ in $H^*$ for $H$ is called cocommutative
if $\Gamma(hh') = \Gamma(h'h)$ for all $h,h' \in H$; if some
integral $\Gamma$ is cocommutative, then all integrals are
cocommutative.

Recall that a Hopf algebra $H$ is called co-Frobenius if $H^{*\rm
rat}$, the unique maximal rational submodule of $H^*$, is nonzero,
or equivalently if  $  \int_l^{H^*} \neq 0$ or if $  \int_r^{H^*}
\neq 0$.

 If $0 \neq \Gamma $ is either a left or a right integral for $H$ in $H^*$,
 then there are bijections from $H$ to $H^{*rat}$ given by
\begin{equation}\label{bij2}
  h\mapsto (h\rh \Gamma) \in H^{*\rm
rat}   \mbox{   and  }     h\mapsto (\Gamma \lh h) \in H^{*\rm
rat}.
 \end{equation}

For $H$ co-Frobenius, a grouplike element in $H$ was defined in
\cite{bdr} as follows.
 Let $0 \neq \lambda
\in \int^{H^*}_l$. If $a\in H$   such that $\l (a)=1$ then by
\cite[Proposition 1.3]{bdr} the element $g^{-1}:= a\leftharpoonup
\lambda = \l (a_1)a_2$ is a grouplike element in $H$ such that for
all $\lambda' \in \int^{H^*}_l$, $\Lambda \in \int^{H^*}_r$, $h \in
H$, we have that
  \begin{equation}\label{g}
\l' (h_1)h_2=\l' (h)g^{-1} \mbox{   and   }  h_1 \Lambda(h_2) =
\Lambda(h) g.
\end{equation}
For the sake of consistency with Radford's definition in
  the finite dimensional case
where the distinguished grouplike in $H$ was constructed using
$\Lambda \in \int_r^{H^*}$, we call $g\in H$ above
  the
distinguished grouplike in $H$.

As well, the following relations hold \cite{bdr}:
\begin{equation}\label{sstarlambda}
\l \circ S=g^{-1}\rh \l~~{\rm ,}~~ \l \circ \smi =\l \lh
g^{-1}~~{\rm and}~~  \l \circ S^2 = g^{-1} \rh \l \lh g.
\end{equation}
The grouplike  $g=1$ if and only if $\int^{H^*}_r = \int^{H^*}_l$.

\subsection{A proof of the $S^4$ formula using integrals}

\vspace{.5cm} There are several short proofs of the formula for
$S^4$ in the literature, for example \cite[p.596]{rad2} gives a
  proof using the trace, Schneider \cite{sch} uses an approach with Frobenius algebras
  and \cite[Theorem 2.10]{moconstanta}
gives a proof combining the Frobenius algebra approach of
Schneider and that of Kadison and Stolin. As well, Radford
\cite{rad2} indicates that a short proof follows from
\cite[Theorem 3 (a),(b)]{rad2}. We supply a  proof for Radford's
$S^4$ formula, based on integrals.

Although the  results in the next lemma are  known (see   \cite{mo}
or \cite{moconstanta} for (i),  \cite[Theorem 3 (a)]{rad2} for an
equivalent form of (ii), and  \cite[Theorem 3 (d)]{rad2} for (iii)),
we supply new proofs for (ii) and (iii).

\begin{lemma}\lelabel{2.1}
Let $H$ be   finite dimensional, $0 \neq t \in \int_l^H$, and $\a $
and $g$ the distinguished grouplike elements of $H^*$ and $H$,
respectively. Let  $0 \neq \Lambda  $ be the unique element of $H^*$
such that  $t \lh \Lambda = \Lambda(t_1)t_2 =1$. Then the following
statements hold:
\begin{itemize}
\item[i)] $\Lambda \in \int^{H^*}_r $; \item[ii)] For any $h\in H$, $h^* \in H^*$ with
 $ h = t \leftharpoonup
h^*  $,  then $ h^* =   \Lambda \lh S(h) =   \a (h_2)\smi (h_1)\rh
\Lambda  $;
\item[iii)] $\Delta (t)=S^2(t_2)g\ot t_1$.
\end{itemize}
\end{lemma}
\begin{proof}
i) The proof of this statement can be found in \cite{mo} or   in
\cite[2.5(i)]{moconstanta}.

ii)Since $t \leftharpoonup \Lambda = 1$, we have $ h(t
\leftharpoonup \Lambda) =h = (t \leftharpoonup \Lambda)h$. Then

$$ h = h\Lambda(t_1)t_2 =\Lambda (S(h_1)h_2t_1)h_3t_2=
\Lambda (S(h)t_1)t_2 =t \lh(\Lambda \lh S(h)).$$

On the other hand, from (\ref{alpha}) we have that $h$ equals
$$
\Lambda (t_1)t_2h=\Lambda (t_1h_2\smi (h_1))t_2h_3 =\a (h_2)\Lambda
(t_1\smi (h_1))t_2 = t \lh (\alpha(h_2)\sinv (h_1) \rh \Lambda).
$$
iii) Let $h^*\in H^*$ with $h^* = \Lambda \lh h$.   By (\ref{g}),
  and $ \Lambda(t_1)t_2 = \Lambda(t)=1$, we have $S(h)g =S(h)\Lambda
  (t)g =S(h)(\Lambda(t_1)t_2)g =
  \Lambda(t_1)S^2(S^{-1}(h)t_2)g$. Now
$$ S(h)\Lambda
  (t)g=S(h)\Lambda (t_2)t_1= \Lambda (h_3t_2)S(h_1)h_2t_1=\le \Lambda \lh h, t_2\ri t_1
  =
  h^*(t_2)t_1 ;
$$
while
$$
\Lambda(t_1)S^2(S^{-1}(h)t_2)g =\Lambda (h_1t_1)S^2(\smi
(h_3)h_2t_2)g=\Lambda (ht_1)S^2(t_2)g=h^*(t_1)S^2(t_2)g.$$

Thus $ h^*(t_2)t_1=h^*(t_1)S^2(t_2)g $, for all $h^*\in H^*$,  and
the statement follows.
\end{proof}

Now we give a short proof of the $S^4$ formula using integrals.

\begin{theorem}\thlabel{2.2}{\rm (\cite[Proposition 6]{rad1})}
Let $H$ be a finite dimensional Hopf algebra with antipode $S$, and $g$ and $\a$
the distinguished grouplike elements of $H$ and $H^*$, respectively. Then
$S^4(h)=g(\a \rh h\lh \a ^{-1})g^{-1}$, for any $h\in H$.
\end{theorem}
\begin{proof}
Let $h\in H$ and let $ h^*\in H^*$ be such that $h=t\lh
h^*=h^*(t_1)t_2$. For $\Lambda   \in \int^{H^*}_r$ as in
\leref{2.1}, using \leref{2.1} ii) and iii), then we have:
 $$
S^2(h)g=h^*(t_1)S^2(t_2)g
 = h^*(t_2)t_1
 =\a (h_2)\Lambda (t_2\smi (h_1))t_1   .
$$

Replacing $h$ with $S^2(h)$, and using $\a \circ S^2=\a $,
(\ref{alpha}), (\ref{g}), and \leref{2.1} i),
  we
obtain
\begin{eqnarray*}
 S^4(h)g &=& (\a \circ S^2)(h_2)\Lambda (t_2S(h_1))t_1 = \a (h_4)\Lambda
(t_2S(h_1))t_1S(h_2)h_3 \\
&=& \a (h_3)\Lambda (t_2)t_1h_2(\a \circ S)(h_1)
 =\a (h_3)gh_2\a ^{-1}(h_1)
=g(\a \rh h\lh \a ^{-1}),
\end{eqnarray*}
and the statement follows. \end{proof}

%%%%%%%%%%%%%%%%%%%%%%%%%%%%%%%%%%%%%%%%%%%%%%%%%%%%%%%%%%%%%%%%
\section{Radford's $S^4$ formula for co-Frobenius Hopf
algebras}\label{cofrob}
%%%%%%%%%%%%%%%%%%%%%%%%%%%%%%%%%%%%%%%%%%%%%%%%%%%%%%%%%%%%%%%%%
\setcounter{equation}{0} Throughout this section, $H$ will denote a
co-Frobenius Hopf algebra, not necessarily finite dimensional. The
following definition will be key in the computations that follow.
 \begin{definition}\delabel{chi} Let $0 \neq \lambda \in \int_l^{H^*}$. Let
  $\chi:\ H\ra H$ be the algebra isomorphism   defined by
\begin{equation}\label{fx1i}
h\rh \l =\l \lh \chi(h),~~\forall~~h\in H.
\end{equation}
 Following the notation for finite
dimensional Frobenius algebras,(see \cite{doitakeuchi}) since
$\lambda(xy) = \lambda(\chi(y)x)$, we call $\chi$ the generalized
Nakayama automorphism for $H$.
\end{definition}

\begin{remark}\relabel{chifd}
For   $H$ finite dimensional,  $\chi(h)= \alpha(h_2)S^{-2}(h_1)$ and
thus  $\epsilon \circ \chi$ is equal to $\alpha$, the distinguished
grouplike in $H^*$. To see this, apply $S^*$ to \leref{2.1}(ii), and
recall  that $\lambda =S^*(\Lambda) \in \int_l^{H^*}$. From
(\ref{sstarlambda}) we obtain $ h \rh \lambda = \alpha(h_2)\lambda
\lh S^{-2}(h_1)$. Finally, apply $\epsilon$ to $\chi(h)=
\alpha(h_2)S^{-2}(h_1)$ to see that $\epsilon \circ \chi = \alpha$.
\end{remark}

 The map $  \epsilon \circ \chi  \in H^*$ is a convolution
invertible algebra map with inverse $ \epsilon \circ \chi \circ
S=\epsilon \circ \chi \circ \smi$. \reref{chifd} justifies the
following terminology.
\begin{definition}\delabel{alpha} For $H$ co-Frobenius and
$\chi$ the generalized Nakayama automorphism, the invertible algebra
map $\epsilon \circ \chi \in H^*$ will be denoted by $\alpha$ and
called the distinguished grouplike in $H^*$.
\end{definition}

\begin{remark}Suppose that $\Lambda$ is a right integral for $H$,
and $\Omega$ is
 the algebra isomorphism defined by $\Omega(h)
\rightharpoonup \Lambda = \Lambda \leftharpoonup h$. Then an easy
computation shows that $\Omega = S^{-1} \circ \chi \circ S = S \circ
\chi \circ S^{-1} $ and so $\epsilon \circ \Omega = \alpha^{-1}$.
\end{remark}

${\;\;\;}$ For $H$  finite dimensional, the space of left integrals
for $H^*$ is invariant under the left adjoint action. This simple
observation motivates the proof of the following generalization of
the formula in \reref{chifd} to co-Frobenius Hopf algebras.

\begin{lemma}\lelabel{mainlemma}
Let $H$ be a co-Frobenius Hopf algebra, $ 0 \neq \l \in
\int^{H^*}_l$
  and $\alpha$ the distinguished grouplike element in
$H^*$.  Then for any $h\in H$,
$$\chi(h)=\a (h_2)S^{-2}(h_1).$$
\end{lemma}
\begin{proof}
We show first that for any $h\in H$ the map $S^2(h_2)\rh \l \lh \smi
(h_1)$ is a left integral. Indeed, for any $h,h'\in H$ we have:
\begin{eqnarray*}
&&\hspace*{-2cm}
\le S^2(h_2)\rh \l \lh \smi (h_1), h'_2\ri h'_1\\
&=&\le \l , \smi (h_1)h'_2S^2(h_2)\ri h'_1\\
&=&\le \l , \smi (h_1)h'_2S^2(h_6)\ri h_3\smi (h_2)h'_1S^2(h_5)S(h_4)\\
\mbox{$(\lambda \in \int^{H^*}_l )$}
&=&\le \l , \smi (h_1)h'S^2(h_4)\ri h_2S(h_3)\\
&=&\le S^2(h_2)\rh \l \lh \smi (h_1), h'\ri 1.
\end{eqnarray*}
Thus for any $h\in H$ there is a scalar $c_h\in k$ such that
$S^2(h_2)\rh \l \lh \smi (h_1)=c_h\l$, and this is equivalent to
\[
\l \lh \chi(S^2(h_2))\smi (h_1)=c_h\l .
\]
By (\ref{bij2}), then  $\chi(S^2(h_2))\smi (h_1)=c_h1$. Applying
$\va$,  we obtain $c_h=\alpha (S^2(h))=\alpha (h)$ so that
$\chi(S^2(h_2))\smi (h_1)=\a (h)1$. Clearly, this is equivalent to
$\chi(S^2(h))=\a (h_2)h_1$ and since $S$ is bijective,
$\chi(h)=\a (h_2)S^{-2}(h_1)$, for any $h\in H$.
\end{proof}

The next corollary generalizes \cite[Corollary 5, p.599]{rad2} to
co-Frobenius Hopf algebras.

\begin{corollary}\colabel{ssqcofrob} For $H$ co-Frobenius with
$0 \neq \lambda \in \int^{H^*}_l$, the following   are equivalent.
\begin{itemize}
\item[i)] $H$ is involutory and $\a = \epsilon$, for $\a$   the
distinguished grouplike in $H^*$. \item[ii)]The integral $\lambda$
is   cocommutative.
\end{itemize}
\end{corollary}
\begin{proof}   Integrals in $\int^{H^*}_l$
are cocommutative if and only if the generalized Nakayama
automorphism $\chi $ is the identity on $H$, so we need only show
that $\chi=id_H$ implies (i).  But for $h \in H$, if $h = \chi(h)=
\a(h_2)S^{-2}(h_1)$,
  then applying $\epsilon$, we obtain that $\epsilon(h) =
  \alpha(h)$ and
  $
  h = S^{-2}(h )$.
\end{proof}

If $H$ is finite dimensional, then applying  \coref{ssqcofrob} to
the dual Hopf algebra $H^*$, we obtain the following statement.

\begin{corollary}\colabel{ssqfinite} Let $H$ be a finite dimensional
Hopf algebra and $0 \neq t
\in \int^H_l$.  Then $t$ is cocommutative if and only if $S^2 =
id_H$ and $g=1$.
\end{corollary}
\begin{proof}Apply   \coref{ssqcofrob} to $H^*$, identifying $H$ with $(H^*)^*$
 and note that
$S^2= id_H$ if and only if $S^{*2}=id_{H^*}$.
\end{proof}

Similar statements to \coref{ssqcofrob} and
 \coref{ssqfinite}  apply for right integrals.

\vspace{.2in}

We   now  prove the main result of this paper.

\begin{theorem}\thlabel{maintheorem}
Let $H$ be a co-Frobenius Hopf algebra, and $\a$, $g$ the
distinguished grouplike elements of $H^*$, $H$, respectively. Then
for any $h\in H$, $$S^4(h)=g(\a \rh h\lh \a ^{-1})g^{-1}.$$
\end{theorem}
\begin{proof}  Let $0 \neq \lambda \in \int_l^{H^*}$. Then by
(\ref{sstarlambda}), $0 \neq \l' = \l \circ S = \l \leftharpoonup
g^{-1} \in \int_l^{(H^{cop)^*}}$. Note that \begin{eqnarray*} h
\rightharpoonup (\lambda \leftharpoonup g^{-1}) &=& (h
\rightharpoonup \l) \leftharpoonup g^{-1}\\
 & =&  \l \leftharpoonup
S^{-2}(h_1)\alpha(h_2)g^{-1} \mbox{ by \leref{mainlemma}}
\\&=& (\l \leftharpoonup g^{-1})
\leftharpoonup (gS^{-2}(h_1)\alpha(h_2)g^{-1}),\end{eqnarray*} so
that the generalized Nakayama automorphism $\chi'$ for $H^{cop}$ is
given by $\chi'(h) = gS^{-2}(h_1)\alpha(h_2)g^{-1}$.
\par Note that $\alpha'$,
the distinguished grouplike for $H^{cop}$ is $\epsilon \circ \chi' =
\alpha$. Also recall that the antipode for $H^{cop}$ is $S^{-1}$.
Now, from \leref{mainlemma} for the co-Frobenius Hopf algebra
$H^{cop}$, we have that $$\chi'(h^{cop}) =
\alpha'(h_2^{cop})S^2(h^{cop}_1) = \alpha(h_1)S^2(h_2).$$ Thus $$
\alpha(h_1)S^2(h_2) =gS^{-2}(h_1)\alpha(h_2)g^{-1}$$ and the
statement follows immediately. \end{proof}

The next corollary yields the analogue of \cite[Theorem
3(b)]{rad2} for the co-Frobenius case. In the finite dimensional
case, this statement can be proved directly and then used together
with \leref{2.1}(ii) or \cite[Theorem 3(a)]{rad2} to give a short
proof of the $S^4$ formula.

\begin{corollary}\colabel{secondchi}
 For any $h \in H$, we have that  $\chi(h)=\a
(h_1)g^{-1}S^2(h_2)g$. \end{corollary}
\begin{proof} From \thref{maintheorem}, $S^2(h) = g(\a \rh S^{-2}(h) \lh
\a^{-1})g^{-1}$, or, equivalently, $S^{-2}(h) = g^{-1}(\a^{-1} \rh
S^{2}(h) \lh \a )g $.  Substituting in the expression for
$\chi(h)$ from \leref{mainlemma} yields the result.
\end{proof}

\begin{corollary}\colabel{finiteorder} If $H$ is co-Frobenius and
 $\alpha$ and $g$ have finite order, then the antipode $S$
of $H$ has finite order.
\end{corollary}

The converse to \coref{finiteorder} fails. In \cite[Example
5.6]{bdgn}, examples were given of infinite-dimensional co-Frobenius
Hopf algebras for which $S^4$ is the identity but the distinguished
grouplike element in $H$ has infinite order. The example below is
based on \cite[Example 5.6(i)]{bdgn} and provides, for any $n$, an
example of an infinite dimensional co-Frobenius Hopf algebra
$\mathcal{H}$ with $S^2$ of order $n$, distinguished grouplike in
$\mathcal{H}$ of infinite order, and distinguished grouplike in
$\mathcal{H}^*$ of order $n$.

\begin{example}\label{bicrossproduct} Let $ T_n$ be the Taft Hopf algebra of dimension
$n^2$.  Then $T_n$ is generated as an algebra by $x,c$ where $c$ is
grouplike and $x$ is $(1,c)$-primitive, i.e., $\triangle(x) = x
\otimes 1 + c \otimes x$. Also $c^n=1$, $x^n = 0$ and $xc = qcx$
where $q$ is a primitive $n$th root of unity.  Let $A = <a>$ be an
infinite cyclic group with identity $e$ and let $kA$ be the group
algebra  . Now, give $T_n$ an $A$-grading by defining $deg(c^ix^j) =
a^j$. Then form the bicrossproduct $\mathcal{H} = T_n \times kA$
with trivial weak action, cocycle and dual cocycle.  The Hopf
algebra $\mathcal{H}$ is equal to the tensor product $T_n \otimes
kA$ as an algebra but has comultiplication defined by
$\triangle(c^ix^j \otimes a^k) = \sum_{0 \leq t \leq j}{j \choose
t}_q (c^ic^t x^{j-t} \otimes a^t a^k)\otimes ( c^ix^t \otimes a^k  )
$. By \cite[Theorem 5.1]{bdgn}, $S(c^ix^j \otimes a^k) =
S_{T_n}(c^ix^j) \otimes S_{kA}(a^{j+k}) $. Thus the order of
$S_{\mathcal{H}}$ is $2n$ since the antipode for $T_n$ has order
$2n$, and the antipode for a group algebra has order 2.  \vspace{.1cm}\\
From the proof of \cite[Proposition 5.2]{bdgn}, $\Lambda =
p_{x^{n-1}} \otimes p_e$ is a nonzero right integral for
$\mathcal{H}$ where $p_h(l) = \delta_{h,l}$.  We compute the order
of the distinguished grouplike in $\mathcal{H}$ using this right
integral. Clearly $(c^ix^j \otimes a^k) \rightharpoonup \Lambda$ is
nonzero only on elements of the form $c^{n-i}x^{n-j-1} \otimes
a^{-k}$.  But we have that $q^i(c^{n-i}x^{n-j-1} \otimes
a^{-k})(c^ix^j \otimes a^k) = (c^ix^j \otimes a^k)(c^{n-i}x^{n-j-1}
\otimes a^{-k})$. Since $q$ is a primitive $n$th root of unity, then
the distinguished grouplike in $\mathcal{H}^*$ has order $n$. By
\cite[Proposition 5.4]{bdgn}, the distinguished grouplike in
$\mathcal{H}$ is $ c^{n-1} \otimes a^{n-1}$, which has infinite
order.\qed
\end{example}

 In general,  for $H$ co-Frobenius, the order of $S$
need not be finite. For example, let $char(k) \neq 2$,  $q \in k$
  not a root of unity, and let $H=k[SL_q(2)]$, the
coordinate ring of $SL_q(2)$ \cite{tak}.  Then $H$ is a cosemisimple
Hopf algebra whose antipode has infinite order. (In fact, by
\cite[Proposition 1.2]{bdgn} no power of the antipode of $H$ is even
inner.) The distinguished grouplike in $H$ is 1.  In the next
example, we compute the map $\chi$ and the distinguished grouplike
element $\a\in H^*$ for $H=k[SL_q(2)]$.
We use the notations of \cite{tak}.\\
${\;\;\;}$ \begin{example}Recall that  $H=k[SL_q(2)]$ is the Hopf
algebra generated as an algebra by $a, b, c, d$ with relations
\[
ba=qab,~~ca=qac,~~db=qbd,~~dc=qcd,~~bc=cb,~~da-qbc=1,~~
ad-q^{-1}bc=1.
\]
The coalgebra structure $\Delta$, $\va$ and the antipode $S$ are
given by
\begin{eqnarray*}
&&\Delta (a)=a\ot a + b\ot c,~~\Delta (b)=a\ot b + b\ot d,\\
&&\Delta (c)=c\ot a + d\ot c,~~\Delta (d)=c\ot b + d\ot d,\\
&&\va (a)= \va(d) =1,~~\va (b)=\va (c)=0, \\
&&S(a)=d,~~S(b)=-qb,~~S(c)=-q^{-1}c,~~S(d)=a.
\end{eqnarray*}
Elements of the form $a^ib^jc^kd^l$ with $i, j, k, l \in
\mathbb{N}$   such that either $i=0$ or $l=0$ form a basis for
$H$.  A left integral    is given by   $\l \in H^*$ defined by
$$ \lambda((bc)^n) = (-1)^n/[n+1] \mbox{ where  } [i] = (q^i -
q^{-i})/(q - q^{-1})$$ and $\lambda $ maps all other basis
elements to 0. Since $H$ is cosemisimple, $g$ = 1.\\
 From the formula for
$S^4$ in \thref{maintheorem}, letting $h=a$, we obtain that
$\alpha(c) = \alpha(b) =0$, and letting $h=b$, we obtain that
$\alpha(d) = q^4 \alpha(a)$. To find $\alpha(a)$, we compute
\begin{eqnarray*}
\le a\rh \l , d\ri
&=&\l (da)=\l (qbc + 1\ri=q\l (bc) + \l (1)\\
&=&-\frac{q(q - q^{-1})}{q^2 - q^{-2}} + 1=\frac{1}{q^2 + 1},\\
\mbox{ and from \leref{mainlemma} we have }\\
\le \l \lh \chi(a), d\ri
&=&\a (a)\l (ad)=\a (a)\l (q^{-1}bc +1)\\
&=&\a (a)\left[-\frac{q^{-1}(q - q^{-1})}{q^2 - q^{-2}} + 1\right]
=\a (a)\frac{q^2}{q^2 + 1}.
\end{eqnarray*}
Thus $\alpha(a) = q^{-2}$. Summarizing, we have that
\begin{eqnarray*}
&&\a (a)=q^{-2},~~\a (b)=\a (c)=0,~~\a (d)=q^2;\\
&&\chi(a)=q^{-2}a,~~\chi(b)=b,~~ \chi(c)=c,~~\chi(d)=q^2d. \qed
\end{eqnarray*}
\end{example}
%%%%%%%%%%%%%%%%%%%%%%%%%%%%%%%%%%%%%%%%%%%%%%%%%%%%%
\section{Integrals and the square of the antipode}
%%%%%%%%%%%%%%%%%%%%%%%%%%%%%%%%%%%%%%%%%%%%%%%%%%%%%
\setcounter{equation}{0} ${\;\;\;}$   Kaplansky's fifth conjecture
 states that a finite dimensional semisimple or cosemisimple Hopf
algebra $H$ is involutory, i.e., that $S^2=id_H$. Over time, this
conjecture has been split into two problems. The first is whether a
semisimple Hopf algebra is cosemisimple and the second  is whether a
semisimple cosemisimple Hopf algebra is involutory.
  The most substantial progress on the first problem was
made by Larson and Radford in \cite{lr1, lr2}. They proved that
over a field of characteristic zero a finite dimensional Hopf
algebra $H$ is semisimple if and only if $H$ is cosemisimple, if
and only if $H$ is involutory. Recently, streamlined proofs using
the Frobenius algebra approach appeared in \cite{moconstanta} for
these results. Larson and Radford also proved that in
characteristic $p$ sufficiently large a semisimple cosemisimple
Hopf algebra is involutory. Ten years later, Etingof and Gelaki
\cite{eg} completed the solution for the second problem. In fact,
using the results of Larson and Radford on the one hand and a
lifting theorem on the other hand, they  proved the following
result.

\begin{theorem}\thlabel{3.0}\cite{eg}
Let $H$ be a finite dimensional Hopf algebra over
 $k$. Then $H$ is
semisimple and cosemisimple if and only if $S^2=id_H$ and $
dim(H)1\not=0$ in $k$.
\end{theorem}

${\;\;\;}$ In this section, we explore equivalent conditions to $H$
involutory  for $H$ co-Frobenius, not necessarily finite
dimensional, in terms of integrals.
\begin{remark}\relabel{remark}
 Let $H$ be
 co-Frobenius with $g=1$. Let $\lambda \in \int^{H^*}_l =
 \int^{H^*}_r$;
  $\lambda
 \circ S = \lambda$. Then the following conditions are equivalent:
\begin{itemize}
\item[(i)] $H$ is involutory.
\item[(ii)] The bilinear form $B(x,y) = \lambda(xS(y))$ is symmetric.
\end{itemize}  For if
 $S^2 =id_H$, then $\lambda(xS(y)) = \lambda \circ S(xS(y)) =
 \lambda(S^2(y)S(x)) = \lambda(yS(x))$. Conversely,
 $\lambda(xS(y)) = \l(yS(x)) = \l \circ S(yS(x)) =\l(S^2(x)S(y))$,
 so that $\l \lh x = \l \lh S^2(x)$ for all $x$, and thus $S^2 =
 id_H$. \end{remark}
In \cite[Corollary 3.6]{larson}, Larson proved that over an
algebraically closed field, a cosemisimple Hopf algebra $H$ is
involutory if and only if $\lambda(h_2 S(h_1)) = \epsilon(h)$ for
some $\lambda \in \int^{H^*}_l$, all $h \in H$.   The proposition
below generalizes this result.

\begin{proposition}\prlabel{larson} For $H$ co-Frobenius over a field $k$,
 the following are equivalent:
\begin{itemize}
\item[i)] $H$ is cosemisimple and involutory. \item[ii)] $H^*$ has
a non-zero left or right integral $\l$ such that
  $\l (S(h_2)h_1)=\va (h)$, for all $h\in
H$. \item[iii)] $H^*$ has a non-zero left or right integral $\l$
such that $\l (h_2S(h_1))=\va (h)$ for all $h\in H$. \item[iv)]
$H$ is cosemisimple and there is a cocommutative integral for $H$
in $H^*$.
\end{itemize}
\end{proposition}
\begin{proof}We show first that (i),(ii) and (iii) are equivalent.
The proof that (i) implies (ii) is the proof that (a) implies (c)
in \cite[Corollary 3.6]{larson}. The implications (i) implies (ii)
and (i) implies (iii)  also follow  immediately from \reref{remark}
and the fact that since $H$ is cosemisimple, a left and right
integral $\lambda$ may be chosen so that $\lambda(1) = 1$.

Now suppose (ii) holds and $\lambda$ is a left or right integral
  with $\l (S(h_2)h_1)=\va (h)$, for all $h\in H$.
 Then $\lambda(1) = 1$ and so $H$
is cosemisimple by the Dual Maschke Theorem and $\lambda  $ is
both a left and right integral. To see that $H$ is involutory, we
compute, for any $h \in H$,
\begin{eqnarray*}
h = h\cdot 1
&=& \lambda(S(h_3)h_2)h_1\\
&=& \lambda(S(h_3)h_2)S^2(h_5)S(h_4)h_1\\
&=& \lambda(S(h_2)h_1)S^2(h_3)\\
&=& S^2(h). \end{eqnarray*}
 Similarly, if (iii) holds, then $H$ is
cosemisimple, and to see that $H$ is involutory, compute $h =
\lambda(h_2S(h_1))h_3 = \lambda(h_4S(h_3))h_5S(h_2)S^2(h_1)$. Thus
each of (ii) and (iii) implies (i). \\
By \coref{ssqcofrob}, it is clear that (iv) implies (i).
Conversely, suppose that the equivalent conditions (i)-(iii) hold.
Note that then the integral in both (ii) and (iii) is the unique
integral $\lambda = \lambda \circ S$ such that $\lambda(1) = 1$.
We show that this implies that $\alpha = \epsilon$. For $h \in H$
and $\lambda = \lambda \circ S$ an integral in $H^*$ such that
$\lambda(1) = 1$, note that
\begin{eqnarray*}
\epsilon(h) = \lambda(h_2S(h_1)) &=& <S(h_1) \rightharpoonup
\lambda
\leftharpoonup h_2,1> \\
&=& <\lambda \leftharpoonup \alpha(S(h_1))S^{-2}(S(h_2))h_3,1>\\
&=& \alpha^{-1}(h_1)<\lambda, S^{-1}(h_2)h_3>\\
&=& \alpha^{-1}(h_1)\lambda(S(h_3)h_2)\\
&=& \alpha^{-1}(h).
\end{eqnarray*}
 Thus $\alpha = \epsilon$ and (iv)
  follows from \coref{ssqcofrob}.\end{proof}

Note that if $H$ is cosemisimple, or more generally if $H$ is
symmetric as a coalgebra \cite{cdn}, then $S^2$ is ``coinner'',
meaning that $S^2(h) = v(h_1)h_2v^{-1}(h_3)$ for some invertible
$v \in H^*$. Then $S^2 = id_H$ if and only if $v*id_H = id_H*v$.

\begin{remark} The equivalent conditions for $S^2 = id_H$ for $H$
co-Frobenius may remind the reader of analogous conditions in the
quasi-triangular and coquasi-triangular case.  It is a well-known
theorem of Drinfel'd \cite{drin} (also proved by Radford
\cite{radQT}) that for $(H,R = R^{(1)} \otimes R^{(2)})$
quasi-triangular, then $S^2$ is an inner automorphism of
$H$, i.e., $S^2(h) = uhu^{-1}$ where $u = S(R^{(2)})R^{(1)}$. \\
Dual results hold for $H$ coquasi-triangular \cite[Lemma
3.3.2]{schauenburg}, \cite[Theorem 1.3]{doi}. Let $(H,\beta)$ be a
coquasi-triangular Hopf algebra and $u \in H^*$ be defined by
$u(h) = \beta(h_2,S(h_1))$. Then $u$ is a unit in $H^*$ and
$S^2(h) = u(h_1)h_2u^{-1}(h_3)$.   These results were generalized
to weakly coquasi-triangular Hopf algebras in \cite{bfm}.
\end{remark}

Now we restrict to the finite dimensional case.

\begin{corollary}\colabel{finite} If $H$ is finite dimensional,
then there exists a non-zero left or
right integral $t$ in $H$ such that $S(t_2)t_1=1$  or such that
$t_2S(t_1) = 1$ if and only if $H$ is semisimple and involutory.
\end{corollary}

${\;\;\;}$ Now, for $H$ finite dimensional, we have a list of
equivalent statements to $H$ involutory.

\begin{theorem}\thlabel{mainss}
Let $H$ be a finite dimensional Hopf algebra over a field $k$ such
that $dim(H)1\not=0$ in $k$. Then the following  are equivalent:
\begin{itemize}
\item[i)] $H$ is semisimple and cosemisimple. \item[ii)] $H$ is
involutory. \item[iii)] $H$ has a non-zero left or right
cocommutative integral $t \in H$. \item[iv)] There exists $ 0 \neq
t \in \int_l^H$ or $0 \neq t \in \int_r^H$   such that
$S(t_2)t_1=1$ or $t_2S(t_1)=1$. \item[v)] $H^*$ has a non-zero
left or right cocommutative integral $\l \in H^*$.
 \item[vi)] $H^*$ has
a non-zero left or right integral $\l $ such that $\l
(h_2S(h_1))=\va (h)$ or $\l (S(h_2)h_1)=\va (h)$, for all $h\in
H$.
\end{itemize}
\end{theorem}
\begin{proof}
The equivalence of i) and ii) follows from \thref{3.0}. By
\coref{ssqfinite}, with the assumption that $\rm{dim}(H)1 \neq 0$,
i) and ii) are equivalent to the existence of a cocommutative
integral in $H$, i.e., to iii). \coref{finite}  implies that i)
and ii) are equivalent to iv). Thus the first four conditions are
equivalent. \\ The equivalence of the remaining conditions is
proved similarly using \coref{ssqcofrob}, \thref{3.0} and
\prref{larson}.
\end{proof}

\begin{remark}
The condition ${\rm dim}(H)1\not=0$ in $k$ in \thref{mainss} is
essential. Let $p$ be   prime, $char(k)=p$   and $C_p=<a>$ a cyclic
group of order $p$. Then $H=kC_p$ is an involutory cosemisimple Hopf
algebra for which ${\rm dim}(H)=0$ in $k$. Moreover, $t=\sum \limits
_{n=0}^{p-1}a^n$ is a non-zero left and right cocommutative integral
for $H$ and $p_e \in H^*$ defined as in Example \ref{bicrossproduct}
is a non-zero left and right cocommutative integral for $H^*$ and
satisfies the conditions in \thref{mainss} vi) by \prref{larson}.
But $H$ is not semisimple and $S(t_2)t_1=t_2S(t_1)=0$. So here ii),
iii), v) and vi) in \thref{mainss} hold, while i) and iv) do not.
Note that in this example, the grouplikes $g$ and $\alpha$ are
trivial, $\lambda(1) = 1$ but $\epsilon(t) = 0$.
\end{remark}

\begin{remark}If   $char(k) = 0$   then by the results of Larson and Radford
mentioned above, $H$ semisimple is equivalent to $H$ cosemisimple
and each of these conditions is equivalent to statements (iii)
through (vi) in   \thref{mainss}.
\end{remark}

${\;\;\;}$ We end this paper by explicitly constructing  the
cocommutative integrals for $H$ and $H^*$ when the finite
dimensional Hopf algebra $H$ is both semisimple and cosemisimple
or, equivalently, when  $H$ is involutory
and ${\rm dim}(H)\not=0$ in $k$.\\
${\;\;\;}$ To this end, consider $\{e_i\}$ a basis in $H$ with dual
basis $\{e^i\}$. We know from \cite{rad2} that the element $r\in H$
defined by
\[
p(r)={\rm Tr}(l(p)\circ S^2),~~\forall ~~p\in H^*
\]
is a non-zero right integral for $H$ if and only if $H$ is cosemisimple,
where we denoted by $l(p)$ the endomorphism of $H$ defined by
$l(p)(h)=p(h_2)h_1$, for any $h\in H$, and by ${\rm Tr}(l(p)\circ S^2)$
the trace of the $H$-endomorphism $l(p)\circ S^2$. So we have
\begin{eqnarray*}
&&\hspace*{-1.5cm}
p(r)={\rm Tr}(l(p)\circ S^2)={\rm Tr}(S^2\circ l(p))\\
&&=\sum \limits _{i=1}^n \le e^i, S^2\left(l(p)(e_i)\right)\ri
=\sum \limits _{i=1}^n \le e^i, S^2((e_i)_1)\ri p((e_i)_2),
\end{eqnarray*}
for all $p\in H^*$, and thus
\[
r=\sum \limits _{i=1}^n\le e^i, S^2((e_i)_1)\ri (e_i)_2.
\]
Replacing $H$ by $H^{\rm op, cop}$ we get that
\[
t=\sum \limits _{i=1}^n \le e^i, S^2((e_i)_2)\ri (e_i)_1
\]
is a non-zero left integral for $H$ if and only if $H$ is cosemisimple. Note that
$\va (r)=\va (t)={\rm Tr}(S^2)$. In particular, if $H$ is also semisimple then $r=t$.\\
${\;\;\;}$
Dually, we have that $\l, \Lambda \in H^*$ defined by
\[
\l =\sum \limits _{i=1}^nS^2(e_i)\rh e^i
~~\mbox{\rm and}~~
\Lambda =\sum \limits _{i=1}^ne^i\lh S^2(e_i)
\]
are non-zero left, respectively right, integrals for $H^*$
if and only if $H$ is semisimple. Clearly, $\l (1)=\Lambda (1)={\rm Tr}(S^2)$, so
if $H$ is also cosemisimple then $\l =\Lambda $.\\
${\;\;\;}$ Now, \thref{mainss} and the uniqueness of integrals for
$H$ and $H^*$ imply:

\begin{corollary}\colabel{3.7}
Let $H$ be a finite dimensional semisimple cosemisimple Hopf
algebra, and $\{e_i\}$ a basis in $H$ with dual basis $\{e^i\}$.
Then
\begin{itemize}
\item[i)] $0 \neq t=\sum \limits _{i=1}^n\le e^i, (e_i)_1\ri
(e_i)_2 =\sum \limits _{i=1}^n \le e^i, (e_i)_2\ri (e_i)_1  $ is a
cocommutative integral; \item[ii)] $0 \neq \l =   \sum \limits
_{i=1}^n e_i\rh e^i=\sum \limits _{i=1}^n e^i\lh e_i  $ is   a
cocommutative integral.
\end{itemize}
\end{corollary}

${\;\;\;}$ \begin{center}  {\bf Acknowledgement:}  Thanks to the
referee for many helpful comments.
\end{center}
%%%%%%%%%%%%%%%%%%%%%%%%%%%%%%%%%%%%%%%%%%%%%%%%%%%%%%%%%%%%

\end{document}